\documentclass[12pt, a4paper]{article}
\usepackage{amssymb, amsmath, amsthm}
\usepackage{multirow}
\usepackage[sort&compress]{natbib}
\usepackage{hyperref}
\usepackage{graphicx}
\usepackage{setspace}
\hypersetup{
    pdftitle={Pooling multiple imputations when the sample happens to be the population},
    pdfauthor={Gerko Vink and Stef van Buuren},
    pdfsubject={Pooling rules for finite populations},
    pdfkeywords={Multiple imputation, pooling, simulation, finite population},
    bookmarksnumbered=true,     
    bookmarksopen=true,         
    bookmarksopenlevel=1,       
    colorlinks=true,
    citecolor=blue,
    linkcolor=blue,
    urlcolor=blue,                 
    pdfstartview=Fit,           
    pdfpagemode=UseOutlines,      
    pdfpagelayout=TwoPageRight
}

\begin{document}
\title{Pooling multiple imputations when the sample happens to be the population.}

\author{Gerko Vink$^{1,2,}$\thanks{Correspondence to: Gerko Vink, Department of Methodology and Statistics, PO Box 80140, 3508 TC Utrecht, the Netherlands. Tel: +31 (0)30 253 9140, Fax: +31 (0)30 253 5797, E-mail: \texttt{G.Vink@uu.nl}}\; and Stef van Buuren$^{1,3}$\\
\newline\\
\small{\textsl{$^1$Department of Methodology and Statistics, Utrecht University}}\\
\small{\textsl{$^2$Division of Methodology and Quality, Statistic Netherlands}}\\
\small{\textsl{$^3$Netherlands Organization for Applied Scientific Research TNO}}}
\date{}
\maketitle

\begin{abstract} 
Current pooling rules for multiply imputed data assume infinite populations. In some situations this assumption is not feasible as every unit in the population has been observed, potentially leading to over-covered population estimates. We simplify the existing pooling rules for situations where the sampling variance is not of interest. We compare these rules to the conventional pooling rules and demonstrate their use in a situation where there is no sampling variance. Using the standard pooling rules in situations where sampling variance should not be considered, leads to overestimation of the variance of the estimates of interest, especially when the amount of missingness is not very large. As a result, populations estimates are over-covered, which may lead to a loss of statistical power. We conclude that the theory of multiple imputation can be extended to the situation where the sample happens to be the population. The simplified pooling rules can be easily implemented to obtain valid inference in cases where we have observed essentially all units and in simulation studies addressing the missingness mechanism only.
\end{abstract}

\section{Background}
Missing data are an ubiquitous problem in medical research. The occurrence of missing data often has an influence on the precision of estimates and may even lead to biased estimates and incorrect statistical inferences. A straightforward approach to obtain valid inference on incomplete data is multiple imputation. With multiple imputation the missing data problem is solved before the analysis takes place and each missing value is imputed $m\geq2$ times, leading to $m$ completed datasets. These $m$ completed datasets are then analyzed separately and their complete data estimates are combined using Rubin's rules \citep{RubinD1987}. 

Current methodology on pooling estimates based on multiply imputed data assumes the data are sampled from infinite populations. In some cases we have data on all units, e.g. rare conditions in medical research and registers in official statistics, and sampling variation plays no role. Yet, even though all units are observed, there may be missing data that affect the precision of the estimates of interest. In such situations, assuming an infinite population may overestimate the variance of the estimates. As a result, confidence intervals are longer than needed, leading to a loss of statistical efficiency. 

This note suggests the use of simplified pooling rules that only account for the variation caused by the mechanism that created the missing data. 

\section*{Methods}
\citet[p. 76]{RubinD1987} defined $Q$ as the quantity of interest (possibly a vector) and $U$ as its variance. With multiple imputation, $m$ complete data estimates can be averaged as
\begin{equation}
\bar{Q}=\frac{1}{m}\sum^{m}_{l=1}{ \hat{Q}}_{l}
\end{equation}
where $\hat Q_l$ is an estimate of $Q$ from the $l$-th imputed
data set. Let $\bar U_l$ be the estimated variance-covariance matrix of
$\hat Q_l$. The complete data variances of $ Q$ can be combined by
\begin{equation}
\bar{U}=\frac{1}{m}\sum^{m}_{l=1}{ {\bar U}}_{l}. 
\end{equation}
The variance between the complete data estimates can be calculated as
\begin{equation}
B=\frac{1}{m-1}\sum^{m}_{l=1}(\hat{ Q}_l-\bar{Q})^\prime(\hat{ Q}_l-\bar{Q}).
\end{equation}
The total variance of $({ Q}-\bar{Q})$ is defined as 
\begin{equation}
T=\bar{U}+B+B/m. 
\end{equation}
For populations for which all units are recorded, the average complete data variance $\bar{U}$ of $ Q$ equals zero - there is no sampling variation - and the total variance of $({ Q}-\bar{Q})$ simplifies to 
\begin{equation}
T=B+B/m. 
\end{equation}
As a consequence, the relative increase in variance due to nonresponse equals
\begin{equation}
r = (1+{m^{-1}}) B/\bar{U}= \infty, 
\end{equation}
and the degrees of freedom $\nu$ can be set to
\begin{equation}
\nu = (m-1)(1+r^{-1})^2=m-1.
\end{equation}

\subsection*{Simulation}
We created a finite population with $N = 1000$ members by drawing 1000 independent realizations from the multivariate normal distribution with means
\begin{equation}\label{meanmatrix}
\mu= \bordermatrix{&	 \cr
 X	&1	\cr
{ Y}_1	&2	\cr
{ Y}_2	&3	},
\end{equation}
and covariance structure
\begin{equation}\label{covmatrix}
\Sigma = \bordermatrix{& X &{ Y}_1 &{ Y}_2	 \cr
 X	&1	&0.1	&0.1\cr
{ Y}_1	&0.1	&1	&0.1\cr
{ Y}_2	&0.1	&0.1	&1},   
\end{equation}
where $ X$ is a completely observed covariate and ${ Y}_1$ and ${ Y}_2$ are made incomplete by randomly deleting values with probabilities that vary between 0.1 and 0.95. 

 Data imputations are performed with \texttt{mice} \cite[version 2.21]{vanbuuren2010} in \texttt{R} \citep[version 3.0.2]{R} with Bayesian linear regression imputation (\texttt{mice.impute.norm}) as the imputation method and 10 iterations for the algorithm to converge. The quantities of scientific interest were the means of $Y_1$ and $Y_2$. The true values were calculated as the sample means before deletion. 
\begin{table}[b!]
\caption{Coverage of the mean. Average results over 10,000 simulations for two variables with varying percentage of missingness. Results are shown for pooling rules for finite populations (simplified rules) and for the pooling rules as defined by Rubin (conventional rules).}
\centering
\begin{tabular}{lrrrrrrrrrrr}
  \hline
 & &\multicolumn{5}{c}{conventional rules}&&\multicolumn{4}{c}{simplified rules}\\ \cline{3-7} \cline{9-12}
 & \%mis & $r$ & $\nu$ & fmi & ciw & cov&& $r$ & $\nu$ & ciw & cov \\ 
  \hline
\multirow{10}{*}{$Y_1$}
   & 10 & 0.13 & 330.36$^*$ & 0.12 & 0.13& 1.000	&& $\infty$ & 4 & 0.06 & 0.949\\  
   & 20 &  0.30 & 775.54 & 0.23 & 0.14 & 1.000    	&&  $\infty$ & 4 &  0.09 & 0.950 \\ 
   & 30 & 0.51 & 401.36 & 0.33 & 0.15 & 0.999  		&&  $\infty$ & 4 &  0.11 & 0.950  \\ 
   & 40 & 0.80 & 400.58 & 0.44 & 0.17 & 0.994  		&&  $\infty$ & 4 &  0.14 & 0.951  \\ 
   & 50 & 1.20 & 57.82 & 0.53 & 0.19 & 0.988     	&&  $\infty$ & 4 &  0.18 & 0.951 \\ 
   & 60 & 1.81 & 31.55 & 0.62 & 0.23 & 0.976    		&&  $\infty$ & 4 &  0.22 & 0.950 \\ 
   & 70 & 2.82 & 20.48 & 0.71 & 0.27 & 0.969    	 	&&  $\infty$ & 4 &  0.27 & 0.953 \\ 
   & 80 & 4.78 & 14.18 & 0.80 & 0.35 & 0.958     	&&  $\infty$ & 4 &  0.35 & 0.952 \\ 
   & 90 & 10.84 & 6.30 & 0.89 & 0.53 & 0.950     	&&  $\infty$ & 4 &  0.54 & 0.951 \\ 
   & 95 & 22.79 & 5.44 & 0.94 & 0.79 & 0.948     	&&  $\infty$ & 4 &  0.80 & 0.951 \\ \cline{2-12}
\multirow{10}{*}{$Y_2$}
   & 10 & 0.13 & 6483.96 & 0.12 & 0.14 & 1.000  	&&  $\infty$ & 4 &  0.06 & 0.948   \\ 
   & 20 & 0.29 & 1028.42 & 0.23 & 0.15 & 1.000 	&&  $\infty$ & 4 &  0.09 & 0.950   \\ 
   & 30 & 0.50 & 261.06 & 0.33 & 0.16 & 0.999 		&&  $\infty$ & 4 &  0.12 & 0.948  	\\ 
   & 40 & 0.78 & 167.69 & 0.43 & 0.18 & 0.995	 	&&  $\infty$ & 4 &  0.15 & 0.949 	\\ 
   & 50 & 1.17 & 90.94 & 0.52 & 0.20 & 0.988 	 	&&  $\infty$ & 4 &  0.18 & 0.951 	\\ 
   & 60 & 1.78 & 68.54 & 0.62 & 0.24 & 0.978 	 	&&  $\infty$ & 4 &  0.23 & 0.949 	\\ 
   & 70 & 2.72 & 19.66 & 0.70 & 0.29 & 0.966 	 	&&  $\infty$ & 4 &  0.28 & 0.949 	\\ 
   & 80 & 4.71 & 11.71 & 0.80 & 0.37 & 0.960 	 	&&  $\infty$ & 4 &  0.37 & 0.951 	\\ 
   & 90 & 10.61 & 8.24 & 0.89 & 0.56 & 0.949 	 	&&  $\infty$ & 4 &  0.57 & 0.949 	\\ 
   & 95 & 22.40 & 5.23 & 0.94 & 0.83 & 0.944 	 	&&  $\infty$ & 4 &  0.84 & 0.947 	\\ 
   \hline
\end{tabular}
\\$^{*}$ Calculated cf. \citet{barnard1999} because occasionally $r\approx0$. \hfill 
\label{tab:results}
\end{table}
\section*{Results}
The results over 10000 simulations are shown in Table \ref{tab:results}. It is clear that excluding the sampling variation in $\bar{U}$ from the confidence interval calculation leads to proper coverage of the 95\% confidence interval of the mean. Taking $\bar{U}$ into account when considering completely observed populations, leads to overcoverage. Not surprisingly, this overcoverage becomes less apparent when the fraction of information missing due to nonresponse approaches 1. 

\begin{figure}[t!]
     		\begin{center}
		\resizebox{\textwidth}{!}{
        			\includegraphics[scale=1]{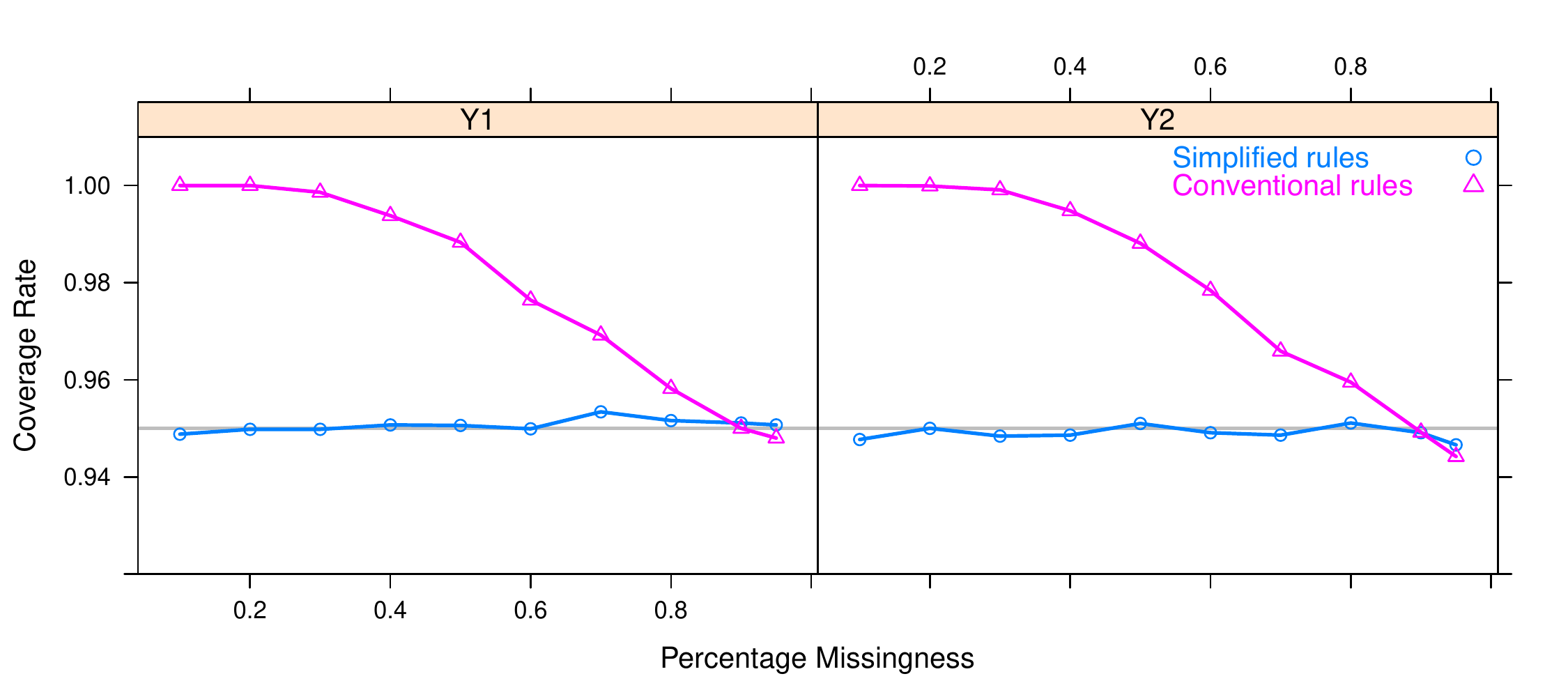}
        		}%
    		\end{center}  
		\caption{Simplified and conventional pooling rules compared. Displayed are coverage rates for different missingness rates when assuming finite (simplified rules) or infinite (conventional rules) populations.}  
  \label{covplot}    
 \end{figure}

With increasing missingness, the role of the between variance $B$ in the total variance $T$ in Rubin's rules becomes increasingly more important and the relative contribution of $\bar{U}$ in $T$ decreases. Due to the increase in $r$, the resulting degrees of freedom $\nu$ approach $m-1$. Eventually, when all data are missing, sampling variation $\bar{U}$ disappears and both pooling approaches become equivalent (see Figure \ref{covplot}). 

All estimates are unbiased. As expected, the conventional pooling rules overestimates the total variance in the datasets covering the entire population, leading to overcoverage (c.f. Table \ref{tab:results}). In contrast, the simplified pooling rules yield consistently a coverage of 95\% of the 95 percent confidence interval.

\section*{Conclusions}
This note illustrates that sharper inferences are possible in situations where the entire population has been observed, and all variation stems for the missing data. Our simulations show that the simplified pooling rules yield variance-covariance estimates that lead to shorter confidence intervals with correct statistical properties. 

Although the adaptation to the conventional pooling rules is small (and may be considered mathematically trivial), we are not aware of any work actually applying simplified pooling. There are several instances where the simplified rules may be of practical interest. First, in situations where essentially all units are observed but missingness has its influence on the precision of estimates, multiple imputation can be utilized to obtain sharper inferences when the proposed pooling rules are used to obtain inference. Such applications can be found throughout many scientific fields, such as medicinal sciences, official statistics and big data applications. 

Another useful application can be found in simulation studies involving the evaluation of imputation approaches. For the last decades, sampling variation has been an essential part of the evaluation of multiple imputation approaches. When infinite populations cannot be assumed during such evaluations, design-based simulation strategies are often used to properly account for sampling variation. However, in order to obtain information about a method's ability to handle the missing data problem, or to objectively compare methods on their ability to correct for missingness, it is not necessary to take sampling variation into account. After all, we are interested only in the missing data mechanism, and are not considering the noise induced by the sampling mechanism for evaluation in such studies. Moreover, not having to consider the sampling mechanism makes the generation of simulation data much more straightforward, especially when generating intricate multivariate data structures.


\bibliographystyle{apalike} 
\bibliography{bibliography}      




%

\end{document}